\input  amstex
\NoBlackBoxes
\documentstyle {amsppt}
\magnification=\magstep1
\baselineskip=24truept
\define \Bxn {B^n_{x_1,x_2,...,x_{k_n}}}
\define \Byn {B^n_{y_1,y_2,...,y_{k_n}}}

\define \Bx {B^n_{x_1,x_2,...,x_k}}
\define \By {B^n_{y_1,y_2,...,y_k}}

\define \subxn{ \{x_1,x_2,...,x_{k_n}\}}
\define \subyn{ \{y_1,y_2,...,y_{k_n}\}}

\define \subx{ \{x_1,x_2,...,x_k\}}
\define \suby{ \{y_1,y_2,...,y_k\}}

\baselineskip=24truept
\NoRunningHeads
\topmatter

\title Law of large numbers for increasing subsequences of random permutations
 \endtitle
\author Ross G. Pinsky\endauthor
\affil Department of Mathematics\\
Technion-Israel Institute of Technology\\
Haifa, 32000
Israel\\
email:pinsky\@math.technion.ac.il
\endaffil
\thanks This research was supported by the Fund for the Promotion
of Research at the Technion and by the V.P.R. Fund.\endthanks

\abstract
Let the random variable $Z_{n,k}$
denote the number
of increasing subsequences of length $k$ in a random
permutation from $S_n$, the symmetric  group of permutations of
$\{1,...,n\}$.
We show that $Var(Z_{n,k_n})=o((EZ_{n,k_n})^2)$
as $ n\to\infty$ if and only if  $k_n=o(n^\frac25)$.
In particular then, the weak law of large numbers holds for
$Z_{n,k_n}$ if  $k_n=o(n^\frac25)$; that is,
$$
\lim_{n\to\infty}\frac{Z_{n,k_n}}
{EZ_{n,k_n}}=1, \ \text{in probability}.
$$
We also show the following approximation
result for the uniform measure $U_n$ on
$S_n$.
Define the probability measure $\mu_{n;k_n}$ on $S_n$ by
$$
\mu_{n;k_n}=\frac1{\binom n{k_n}}\sum_{x_1<x_2<...<x_{k_n}} U_{n;x_1,x_2,...x_{k_n}},
$$
where
$U_{n;x_1,x_2,...,x_{k_n}}$ denotes the uniform measure
on the subset of permutations which contain the increasing
subsequence $\{x_1,x_2,...,x_{k_n}\}$.
Then the weak law of large numbers holds for $Z_{n,k_n}$
if and only if
$$
\lim_{n\to\infty}||\mu_{n;k_n}-U_n||=0,\tag*
$$
where $||\cdot||$ denotes the total variation norm.
In particular then, (*) holds if $k_n=o(n^\frac25)$.

In order to evaluate the asymptotic behavior of the second moment,
we need to analyze  occupation times of certain conditioned
two-dimensional random walks.

\endabstract
\endtopmatter
\pagebreak
\noindent \bf 1. Introduction and Statement of Results.\rm\
Let $S_n$ denote the symmetric group of permutations
of  $\{1,...,n\}$.
By introducing the uniform probability measure $U_n$ on $S_n$,
one can consider $\sigma\in S_n$ as a random permutation.
Probabilities and expectations according to $U_n$ will frequently
be denoted by
the generic notation $P$ and $E$ respectively.
The problem of analyzing  the distribution of the length, $L_n$, of the longest
increasing subsequence  in a random permutation
from $S_n$
has a long and
distinguished history; see \cite{1} and references therein. In
particular, the work of  Logan and Shepp
\cite{6} together with that   of Vershik and Kerov \cite{8}
show that $EL_n\sim 2n^\frac12$ and that
$\sigma^2(L_n)=o(n)$,
as $n\to\infty$.
Profound recent work by Baik, Deift and Johansson \cite{2}
has shown that $\lim_{n\to\infty}P(\frac{L_n-2n^\frac12}{n^\frac16}\le x)=
F(x)$, where $F$ is an explicitly identifiable function.

There doesn't seem to be any literature on the
random variable $Z_{n,k}=Z_{n,k}(\sigma)$,
which we define to be the number
of increasing subsequences of length $k$ in a permutation $\sigma\in S_n$.
Thus, for example, if $\sigma=\left(\matrix
1&2&3&4&5\\1&3&4&5&2\endmatrix\right)$,
then $Z_{5,3}(\sigma)=4$ because there are four increasing subsequences of length
three; namely, 134, 135, 145 and 345.
It is useful to represent $Z_{n,k}$ as a sum of indicator random variables.
For positive integers $\{x_1,...,x_k\}$ satisfying
$1\le x_1<x_2<...<x_k\le n$,  let $B^n_{x_1,...x_k}
\subset S_n$ denote the subset of permutations which contain the increasing
subsequence $\{x_1,x_2,...,x_k\}$. Then we have
$$
Z_{n,k}=\sum_{x_1<x_2<...<x_k}1_{B^n_{x_1,x_2,...,x_k}},
$$
where the sum is over the $\binom nk$ distinct increasing
subsequences of length $k$. Since the probability that a random
permutation fixes any particular increasing sequence of length $k$
is $\frac1{k!}$, it follows that the expected value of $Z_{n,k}$
is given  by
$$
EZ_{n,k}=\frac{\binom nk}{k!}.\tag1.1
$$

One can consider $k$ to depend on $n$ in which case we write $k_n$.
We are interested in a law of large numbers of the form
$\frac{Z_{n,k_n}}
{EZ_{n,k_n}}\rightarrow 1$ in probability, for
appropriate choices of $k_n$.
Of course, in light of the above cited works on the longest increasing subsequence,
such a result cannot hold for $k_n\ge cn^\frac12$ with $c>2$.
A straightforward calculation using Stirling's formula shows that
$$
\aligned&EZ_{n,cn^l}\sim\frac1{2\pi cn^l}\left[(\frac ec)^2n^{1-2l}\right]^{cn^l},\ \text{as}\ n\to\infty,
\ \text{for}\ l\in(0,\frac12);\\
&EZ_{n,cn^\frac12}\sim\frac{\exp(-\frac{c^2}2)}{2\pi cn^\frac12}(\frac ec)^{2cn^\frac12}.
\endaligned\tag1.2
$$
(For the case $k_n=cn^\frac12$, we have used the fact that
$\lim_{n\to\infty}\prod_{j=0}^{cn^\frac12-1}(1-\frac jn)=\exp(-\frac{c^2}2)$,
which is proved by taking the logarithm of the above product.
Note that the factor $\exp(-\frac{c^2}2)$ suddenly appears in the formula when $l=\frac12$.)
In particular then, it follows from (1.2) that
$\lim_{n\to\infty}EZ_{n,k_n}=\infty$, if $k_n\le cn^\frac12$ with $c<e$,
and $\lim_{n\to\infty}EZ_{n,k_n}=0$, if $k_n\ge en^\frac12$.

The law of large numbers for $Z_{n,k_n}$ is in fact  equivalent
to a certain approximation result for the uniform measure, which
we now describe.
Recall that for probability measures $P_1$ and $P_2$ on $S_n$, the total
variation norm is defined by
$$
||P_1-P_2||\equiv\max_{A\subset S_n}(P_1(A)-P_2(A))=\frac12\sum_{\sigma\in S_n}|P_1(\sigma)-P_2(\sigma)|.
$$
For $x_1<x_2<...<x_{k_n}$, let $U_{n;x_1,x_2,...,x_{k_n}}$
denote the uniform measure on  permutations
which have $\{x_1,x_2,...,x_{k_n}\}$ as an increasing sequence; that
is $U_{n;x_1,x_2,...,x_{k_n}}$ is uniform on $B^n_{x_1,x_2,...,x_{k_n}}$.
Note that $U_{n;x_1,x_2,...,x_{k_n}}$ is defined by
$U_{n;x_1,x_2,...,x_{k_n}}(\sigma)=\frac{k_n!}{n!}1_{B^n_{x_1,x_2,...,x_{k_n}}}(\sigma)$.
Now define the probability measure $\mu_{n;k_n}$ on $S_n$ by
$$
\mu_{n;k_n}=\frac1{\binom n{k_n}}\sum_{x_1<x_2<...<x_{k_n}} U_{n;x_1,x_2,...x_{k_n}}.
$$
Equivalently,
$$
\mu_{n;k_n}(\sigma)=\frac1{\binom n{k_n}}\frac{k_n!}{n!}Z_{n,k_n}(\sigma),\
\ \sigma\in S_n.\tag1.3
$$
The measure $\mu_{n;k_n}$ can be realized concretely as follows.
Consider $n$ cards, numbered from 1 to $n$, and laid out on a table from left
to right in increasing order. Place a black mark on $k_n$ of the cards,
chosen at random. Pick up all the cards without black marks and then
randomly insert them between the $k_n$ cards with black marks that
remain on the table.
The resulting distribution is $\mu_{n;k_n}$.

\medskip

\bf\noindent Proposition 1.\it\ The  law of large numbers holds for
$Z_{n,k_n}$; that is
$$
\lim_{n\to\infty}\frac{Z_{n,k_n}}{EZ_{n,k_n}}=1\
\text{in probability},
$$
if and only if
$$
\lim_{n\to\infty}||\mu_{n;k_n}-U_n||=0.\tag1.4
$$
\medskip\rm

The proof of Proposition 1 appears at the end of this section.

The measure  $\mu_{n;k_n}$ corresponds to ignoring a set of $k_n$
random cards and randomizing the rest of the cards.
How many random  cards can one afford to ignore like this
and maintain asymptotic randomness? Corollary 2 below
shows that one can afford to ignore $k_n=o(n^\frac25)$ cards,
while  the results cited  above on the longest increasing subsequence
show that one certainly cannot  afford to ignore $cn^\frac12$ cards
for $c>2$.

For the law of large numbers we will use Chebyshev's inequality.
The calculation of the  second moment  is  nontrivial because it
involves expectations of the form $E1_{B^n_{x_1,...,x_{k_n}}}
1_{B^n_{y_1,...,y_{k_n}}}$, and these expectations depend rather
intimately on the relative positions of $\{x_1,x_2,....,x_{k_n}\}$
and $\{y_1,y_2,....,y_{k_n}\}$.
We begin with the explicit form of the second moment of
$Z_{n,k}$ for any $k\le n$.

\medskip

\noindent \bf Proposition 2.\it
$$
EZ^2_{n,k}=\sum_{j=0}^k\binom n{2k-j}\frac1{(2k-j)!}
A(k-j,j),
$$
where
$$
A(N,j)=\sum\Sb \sum_{r=0}^jl_r=N\\ \sum_{r=0}^jm_r=N\endSb
\prod_{r=0}^j\left(\frac{(l_r+m_r)!}{l_r!m_r!}\right)^2.\tag1.5
$$

\medskip\rm

In order to evaluate the asymptotic behavior of Var($Z_{n,k_n}$), one
must be able to adequately evaluate the asymptotic behavior of
$A(k_n-j,j)$.
In fact, it turns out that we need a good lower bound
for $A(k_n-1,1)$ and a good upper bound for $A(k_n-j,j)$, for all $j=1,2,...k_n$.
We were able to interpret  $\frac{A(N,j)}{\binom{2N}N^2}$
as the sum of  certain expected occupation times of the horizontal axis
for the standard, simple, symmetric two-dimensional  random walk
starting from the origin and conditioned on returning to the
origin  at the $2N$-th step.
This characterization was sufficient to obtain the appropriate
bounds to prove the following theorem.

\medskip

\noindent \bf Theorem 1.\it\

\noindent i. If $k_n=o(n^\frac25)$, then
$$
\frac{Var(Z_{n,k_n})}{(EZ_{n,k_n})^2}=O(\frac{k_n^\frac52}n),
\ \text{as}\ n\to\infty;
$$
In particular then,
$Var(Z_{n,k_n})=o((EZ_{n,k_n})^2), \ \text{as}\ n\to\infty$.

\noindent ii. If   $c_1n^{\frac25}\le k_n\le c_2n^\frac25$,
for constants $c_1,c_2>0$, then
$$
c_3(EZ_{n,k_n})^2\le Var(Z_{n,k_n})\le c_4(EZ_{n,k_n})^2,
$$
for constants $c_3,c_4>0$.

\noindent iii. If $\lim_{n\to\infty}n^{-\frac25}k_n=\infty$
and $\limsup_{n\to\infty}n^{-\frac12}k_n<\infty$, then
$$
\lim_{n\to\infty}\frac{Var(Z_{n,k_n})}{(EZ_{n,k_n})^2}=\infty.
$$

\medskip\rm

\noindent \bf Corollary 1.\it\ i. If $k_n=o(n^\frac25)$,
then
$$
\lim_{n\to\infty}\frac{Z_{n,k_n}}
{EZ_{n,k_n}}=1, \ \text{in probability};
$$

\noindent ii. If $k_n=O(n^\frac25)$, then
$$
\liminf_{n\to\infty}P(\frac{Z_{n,k_n}}{EZ_{n,k_n}}>\delta)>0,
\ \text{for some}\ \delta>0.
$$
\medskip\rm
Part (i) of Corollary 1 follows immediately from Chebyshev's inequality and Theorem 1-i.
The proof of part (ii) of Corollary 1 appears below.

Corollary 1 and Proposition 1 yield immediately  the following  approximation result.
\medskip

\noindent \bf Corollary 2.\it\ If $k_n=o(n^\frac25)$, then
$$
\lim_{n\to\infty}||\mu_{n;k_n}-U_n||=0.
$$

\medskip\rm
In light of the above results, we pose the following question:

\bf \noindent Open Question:\it\ Presumably there exists a critical
exponent $l_c$ such that
the law of large numbers holds
for $Z_{n,n^l}$ with $l<l_c$ and does not hold for $l>l_c$.
What is $l_c$?
\rm\

In section two we prove Proposition 2 and in section 3 we prove
Theorem 1. Lemmas 2 and 3, which appear in section 3 and
 give the key estimates on $A(N,j)$
used in  the proof of Theorem 1, are proved in section four.

The literature on  increasing subsequences in random
permutations in a context other than that of the largest such
subsequence is very scarce. The
 random variable $Z_n$, defined as the total number of increasing subsequences
of all possible lengths in a random permutation,
was studied in \cite{5}. Both $EZ_n$ and $Var(Z_n)$ were calculated explicitly and evaluated
asymptotically. It turns out that $Var(Z_n)$ is of a larger
order than $(EZ_n)^2$, so it is not possible to apply Chebyshev's
inequality and obtain a law of large numbers.
However, the authors were able to show that
$\frac{\log Z_n}{n^\frac12}$ converges in probability and in
mean to a positive constant.
In \cite{3}, the random variable $Z_{n,k}$ actually
appears in a different guise. Equation (1.1) appears
there as well as an upper bound for $EZ_{n,cn^\frac12}$; however, this
random variable  is not the object of study in that paper. In \cite{7},
\it inversions\rm\---which are decreasing subsequences of length 2---are studied, and a central limit
theorem is proved.

We conclude this section with the proofs of
Corollary 1-ii and  Proposition 1.
\medskip

\noindent \bf Proof of Corollary 1-ii. \rm\
Assume to the contrary that the result is not true. Then there
exists a subsequence  $\{(n_i,k_{n_i})\}_{i=1}^\infty$ of
$\{(n,k_n)\}_{n=1}^\infty$,
such that
$\frac{Z_{n_i,k_{n_i}}}{EZ_{n_i,k_{n_i}}}$ goes to 0 in probability.
By taking a further subsequence if necessary, we may assume that
either $k_{n_i}=o(n_i^\frac25)$ as $i\to\infty$, or
$\lim_{i\to\infty}n_i^{-\frac25}k_{n_i}=c>0$.
In light of part (i), we
obtain a contradiction in the former case. Thus,
it remains to consider the latter case.
In this case, it follows from Theorem 1-ii that $\frac{(EZ_{n_i,k_{n_i}})^2}
{Var(Z_{n_i,k_{n_i}})}$
is bounded away from 0 and $\infty$. Using this along with the assumption that
$\frac{Z_{n_i,k_{n_i}}}{EZ_{n_i,k_{n_i}}}\rightarrow 0$ in probability, we conclude that
$$
\lim_{n\to\infty}P(\frac{Z_{n_i,k_{n_i}}-EZ_{n_i,k_{n_i}}}{\sqrt {Var(Z_{n_i,k_{n_i}})}}\le-\rho)=1,\tag1.6
$$
for some $\rho>0$.
However, since
the second moments of the
$\frac{Z_{n_i,k_{n_i}}-EZ_{n_i,k_{n_i}}}{\sqrt {Var(Z_{n_i,k_{n_i}})}}$
are equal to 1, this  quotient is uniformly bounded. The uniform boundedness along with (1.6)
contradict  the fact that the first moment of
$\frac{Z_{n_i,k_{n_i}}-EZ_{n_i,k_{n_i}}}{\sqrt{ Var(Z_{n_i,k_{n_i}})}}$ is 0.
\hfill $\square$
\medskip

\noindent \bf Proof of Proposition 1.\rm\
For $\epsilon\in(0,1)$, define
$$
D_{n,\epsilon,k_n}=\{\sigma\in
S_n:\frac{\mu_{n;k_n}(\sigma)}{U_n(\sigma)}\in[1-\epsilon,1+\epsilon]\}.
$$
We claim that
(1.4) holds if and only if
$$
\lim_{n\to\infty}U_n(D^c_{n,\epsilon,k_n})=0,\ \text{for all}\ \epsilon>0.\tag1.7
$$
We first show the sufficiency of (1.7). Since
$\lim_{n\to\infty}U_n(D_{n,\epsilon,k_n})=1$,
 it follows from
the definition of $D_{n,\epsilon,k_n}$ that
$\liminf_{n\to\infty}\mu_{n;k_n}(D_{n,\epsilon,k_n})\ge1-\epsilon$, and thus
$$
\limsup_{n\to\infty}\mu_{n;k_n}(D^c_{n,\epsilon,k_n})\le \epsilon.\tag1.8
$$
Thus, for any $A_n\subset S_n$, we have
$$
\aligned&|U_n(A_n)-\mu_{n;k_n}(A_n)|\le |U_n(A_n\cap D_{n,\epsilon,k_n})-\mu_{n;k_n}(A_n\cap D_{n,\epsilon,k_n})|\\
&+|U_n(A_n\cap D^c_{n,\epsilon,k_n})-\mu_{n;k_n}(A_n\cap D^c_{n,\epsilon,k_n})|.
\endaligned\tag1.9
$$
By the definition of $D_{n,\epsilon,k_n}$, the first term on the
right hand side of (1.9) is no greater than $\epsilon$. By (1.7)
and (1.8), the $\limsup$ of the second term on the right hand side of (1.9)
is no greater than $\epsilon$.
Since $\epsilon>0$ is arbitrary, we conclude that
$$
\lim_{n\to\infty}|U_n(A_n)-\mu_{n;k_n}(A_n)|=0.
$$
Since the sets $\{A_n\}$ are arbitrary, this
proves (1.4).

We now show the necessity of (1.7). Let
$$
C_{n,\epsilon,k_n}=
\{\sigma\in
S_n:\frac{\mu_{n;k_n}(\sigma)}{U_n(\sigma)}<1-\epsilon\}.
$$
If (1.7) does not hold, then we may assume without loss of generality that
there exists a $\delta>0$ and an $\epsilon_0>0$ such that
$U_n(C_{n,\epsilon_0,k_n})\ge \delta$, for all $n$.
But then, from the definition of $C_{n,\epsilon_0,k_n}$, it follows that
$\mu_{n;k_n}(C_{n,\epsilon_0,k_n})<(1-\epsilon_0)U_n(C_{n,\epsilon_0,k_n})$,
and thus $|U_n(C_{n,\epsilon_0,k_n})-\mu_{n;k_n}(C_{n,\epsilon_0,k_n})|>\epsilon_0\delta$,
for all $n$, which shows that (1.4) does not hold.

To complete the proof of the proposition then, it remains to prove that (1.7)
holds if and only if the law of large numbers holds.
Using (1.3) for the first equality below, and using
(1.1) for the second equality, we
have
$$
\aligned&U_n(D_{n,\epsilon,k_n}^c)=
P(\frac1{\binom n{k_n}}\frac{k_n!}{n!}Z_{n,k_n}
\not\in[\frac{1-\epsilon}{n!},\frac{1+\epsilon}{n!}])=\\
&=P(|\frac{Z_{n,k_n}}{EZ_{n,k_n}}-1|>\epsilon).
\endaligned\tag1.10
$$
From (1.10), it follows that
(1.7) holds if and only if the law of large numbers holds for
$Z_{n,k_n}$.
\hfill $\square$

\medskip

\noindent \bf 2. Proof of Proposition 2.\rm\
From the definition of $Z_{n,k}$, it follows that
$$
EZ^2_{n,k}=\sum E1_{\Bx}1_{\By},\tag2.1
$$
where the sum is over the $\binom n{k}^2$ pairs
$\Bx,\By$ with $x_1<x_2<...<x_{k}$ and
$y_1<y_2<...<y_{k}$.
It turns out that
$E1_{\Bx}1_{\By}$ depends rather intimately on the relative positions
of $\subx$ and $\suby$.

Let $j\in\{0,1,...,k\}$.
For any particular subset $A\subset\{1,2,..., n\}$ satisfying $|A|=2k-j$,
there are   $\binom {2k-j}j\binom{2k-2j}{k-j}$ ordered pairs of sets
$\Bx$, $\By$ for which $\{x_1,x_2,...,x_{k}\}\cup\{y_1,y_2,...,y_{k}\}
=A$. Of course, it follows that $|\{x_1,x_2,...,x_{k}\}
\cap\{y_1,y_2,...,y_{k}\}|=j$.
We will say that such a  pair $\Bx,\By$ \it corresponds to\rm\ $A$.
For any pair $\Bx,\By$ corresponding to $A$, there exist
numbers $\{l_r\}_{r=0}^j$ and $\{m_r\}_{r=0}^j$  such that
exactly $l_0$ elements of $\subx$ and  $m_0$
elements of $\suby$ strictly precede the first element that is common to  $\subx$
and $\suby$, exactly $l_r$ elements of $\subx$ and $m_r$ elements of $\suby$
fall strictly between the $r$-th and the $(r+1)$-th element that is common to
$\subx$ and $\suby$, for $r=1,2,...,j-1$, and exactly $l_j$ elements
of $\subx$ and $m_j$ elements of $\suby$ strictly follow the $j$-th and final
element that is common to $\subx$ and $\suby$.
We will refer to the
numbers $\{l_r\}_{r=0}^j$ and
$\{m_r\}_{r=0}^j$
as the ``interlacing numbers'' for the pair  $\Bx,\By$.

\medskip

\noindent \bf Lemma 1.\it \ Let the pair $\Bx,\By$ satisfy
$$
|\subx\cap\suby|=j
$$
and let $\{l_r\}_{r=0}^j$, $\{m_r\}_{r=0}^j$ be the corresponding
interlacing numbers. Then
$$
E1_{\Bx} 1_{\By}
=\frac1{(2k-j)!}\prod_{r=0}^j\frac{(l_r+m_r)!}{l_r!m_r!}.\tag2.2
$$

\medskip\rm

\noindent \bf Proof.\rm\
Without loss of generality, we may assume that $n=2k-j$, since
only the relative positions of the $2k-j$
distinct points in the set
$\{x_1,x_2,...,x_{k},y_1,y_2,...,y_{k}\}$
are relevant.
Thus, we are considering permutations from $S_{2k-j}$.
For each $r=0,1,...,j$, consider the $l_r+m_r$
positions between common elements (of course, for $r=0$ and $r=j$, ``between''
is not the correct word). There are $(l_r+m_r)!$ ways to fill these positions.
However, if we require that the $l_r$ positions reserved for the $x$-chain
and the $m_r$ positions reserved for the $y$-chain  be in increasing order, this
reduces the number of ways to $\frac{(l_r+m_r)!}{l_r!m_r!}$.
Thus, there are
$\prod_{r=0}^j\frac{(l_r+m_r)!}{l_r!m_r!}$ ways to fill all the positions so that
$\Bx\cap \By$ will occur, and of course, all together there are $(2k-j)!$ ways to fill
the positions with no restrictions.\hfill $\square$

We now complete the proof of the proposition.
Simple combinatorial considerations show  that
out of the $\binom {2k-j}j\binom{2k-2j}{k-j}$ pairs
$\Bx$, $\By$ corresponding to a set $A$ satisfying  $|A|=2k-j$,
there are $\prod_{r=0}^j\binom{l_r+m_r}{l_r}=\prod_{r=0}^j\frac
{(l_r+m_r)!}{l_r!m_r!}$
of them with the interlacing numbers $\{l_r\}_{r=0}^j$, $\{m_r\}_{r=0}^j$.
Using this fact along with (2.1), (2.2)
and the fact that
there are $\binom n{2k-j}$ distinct subsets $A\subset\{1,2,...,n\}$
such that $|A|=2k-j$,
we obtain the
formula for $EZ^2_{n,k}$ in Proposition 2.
\hfill $\square$

\medskip

\noindent \bf 3. Proof of Theorem 1.\rm\
Similar to (2.1),
we can write the variance of
$Z_{n,k_n}$ in the form
$$
Var(Z_{n,k_n})=\sum E1_{\Bxn}1_{\Byn}-\sum E1_{\Bxn}E1_{\Byn},\tag3.1
$$
where the sum is over the $\binom n{k_n}^2$ pairs
$\Bxn,\Byn$.
By (2.2),
$$
\aligned&E1_{\Bxn}1_{\Byn}- E1_{\Bxn}E1_{\Byn}=0, \ \text{if}\\
&\subxn \ \text{and}\ \subyn \ \text{are disjoint}.\endaligned\tag3.2
$$
The number of pairs $\subxn$, $\subyn$ which are not disjoint is equal
to $\binom n{k_n}^2-\binom n{k_n}\binom{n-k_n}{k_n}$.
If $k_n=o(n^\frac12)$, then a simple calculation
reveals that \linebreak
$\binom n{k_n}^2-\binom n{k_n}\binom{n-k_n}{k_n}=o(\binom n{k_n}^2)$.
Thus,
$$
\aligned&\sum_{\subxn\cap\subyn\neq\emptyset} E1_{\Bxn}E1_{\Byn}\\
&=\left(\binom n{k_n}^2-\binom n{k_n}
\binom {n-k_n}{k_n}\right)\frac1{(k_n!)^2}=
o(\frac{\binom n{k_n}^2}{(k_n!)^2})=
o((EZ_{n,k_n})^2),
\endaligned\tag3.3
$$
where the final equality follows from (1.1).
On the other hand,
if it is not true that $k_n=o(n^\frac12)$, then the left hand side
of (3.3) will be $O((EZ_{n,k_n})^2)$.
In light of this last remark along with (3.1)-(3.3), the
theorem will be proved once we    show that
$$
\aligned&\frac{\sum_{\subxn\cap\subyn\neq\emptyset} E1_{\Bxn}1_{\Byn}}
{(EZ_{n,k_n})^2}
=O(\frac{k_n^\frac52}n),\\
& \ \text{if}\ k_n\ \text{is as in part (i)};\endaligned
\tag3.4-a
$$
$$
\aligned&\frac{\sum_{\subxn\cap\subyn\neq\emptyset} E1_{\Bxn}1_{\Byn}}
{(EZ_{n,k_n})^2}\ \\
&\text{is bounded from 0 and }\ \infty\ \text{ if}\
k_n\ \text{is as in part (ii)};\endaligned\tag3.4-b
$$
$$
\aligned&\lim_{n\to\infty}\frac{\sum_{\subxn\cap\subyn\neq\emptyset} E1_{\Bxn}1_{\Byn}}
{(EZ_{n,k_n})^2}=\infty,\\
&  \text{if}\ k_n\ \text{is as in part (iii)}.
\endaligned\tag3.4-c
$$

By Proposition 2 and its proof, it follows that
$$
\aligned&\sum_{\subxn\cap\subyn\neq\emptyset} E1_{\Bxn}1_{\Byn}\\
&=\sum_{j=1}^{k_n}\binom n{2k_n-j}\frac1{(2k_n-j)!}
A(k_n-j,j),\ \text{where} \ A(N,j) \ \text{is as in (1.5)}.\endaligned
$$
Using this with (3.4)
and the fact that
$EZ_{n,k_n}=\frac{\binom n{k_n}}{k_n!}$, the proof will be complete if we show
that
$$
\frac{(k_n!)^2}{\binom n{k_n}^2}
\sum_{j=1}^{k_n}\binom n{2k_n-j}\frac1{(2k_n-j)!}A(k_n-j,j)=O(\frac{k_n^\frac52}n),
\ \text{if }\ k_n\ \text{is as in part (i)};\tag3.5-a
$$
$$
\aligned&\frac{(k_n!)^2}{\binom n{k_n}^2}
\sum_{j=1}^{k_n}\binom n{2k_n-j}\frac1{(2k_n-j)!}A(k_n-j,j)
\ \text{is bounded from}\ 0 \ \text{and}\ \infty \\
&\text{if}\
k_n\ \text{is as in part (ii)}\ ;\endaligned\tag3.5-b
$$
$$
\lim_{n\to\infty}\frac{(k_n!)^2}{\binom n{k_n}^2}
\sum_{j=1}^{k_n}\binom n{2k_n-j}\frac1{(2k_n-j)!}A(k_n-j,j)
=\infty,\ \text{if}\ k_n\ \text{is as in part (iii)}.
\tag3.5-c
$$

It remains therefore to analyze the left hand side of (3.5).
In the next section we will prove the following key estimates:

\medskip

\noindent \bf Lemma 2.\it\ For each
$\rho\in(0,\infty)$, there exists a constant $C_\rho>0$
such that
$$
A(N,j)\le
C_\rho^j\ \frac {j^\frac12}
{\Gamma(\frac {j+1}2)}(2N)^\frac j2
\binom{2N}N^2,\ \text{for}\ j,  N\ge1 \
\text{and}\ \frac jN\le\rho.
$$
In particular, since $A(N,j)$ is increasing in $N$, one has
$$
A(k-j,j)\le A(k,j)\le C_1^j\frac{j^\frac12}{\Gamma(\frac {j+1}2)}
(2k)^\frac j2\binom{2k}k^2, \ \text{for}\  j,  k\ge1 \ \text{and}\ j\le k.
$$
\medskip
\rm

\noindent\bf Lemma 3.\it\
There exists a constant $C>0$ such that
$$
A(k-1,1)\ge C(2k-2)^\frac12\binom {2k-2}{k-1}^2.
$$
\medskip\rm

We now use Lemma 2 to show that (3.5-a) and the part of (3.5-b)
concerning boundedness  from $\infty$ hold.
Afterwards, we will use Lemma 3 to show  that (3.5-c) and
the part of (3.5-b)
concerning boundedness  from 0 hold.

In light of Lemma 2, it suffices to  show that (3.5-a) and
the part of (3.5-b)
concerning boundedness from $\infty$
hold with
$A(k_n-j,j)$ replaced by
$C_1^j\frac{j^\frac12}{\Gamma(\frac {j+1}2)}(2k_n)^\frac j2
\binom{2k_n}{k_n}^2$.
Letting
$$
B(n,k_n,j)=
\frac{(k_n!)^2}{\binom n{k_n}^2}\binom n{2k_n-j}\frac1{(2k_n-j)!}
\binom{2k_n}{k_n}^2
(2k_n)^\frac j2,
$$
it follows that (3.5-a) (respectively the part of (3.5-b) concerning
boundedness from $\infty$) will hold if we show that
$$
\sum_{j=1}^{k_n}B(n,k_n,j)\frac{C_1^jj^\frac12}{\Gamma(\frac {j+1}2)}
$$
is $O(\frac{k_n^\frac52}n)$, if $k_n$ is as in part (i) (respectively, bounded
if $k_n$ is as in part (ii)).
Simplifying and making some cancellations,
we have
$$
B(n,k_n,j)=\frac{((n-k_n)!)^2}{n!(n-2k_n+j)!}
\left(\frac{(2k_n)!}{(2k_n-j)!}\right)^2(2k_n)^\frac j2.\tag3.6
$$
We have
$$
b_1n^{-j}\le \frac{((n-k_n)!)^2}{n!(n-2k_n+j)!}\le b_2n^{-j},
\ j=1,...,k_n,
\tag3.7
$$
for positive constants $b_1,b_2$.
(For the lower bound, we have used
the fact that $k_n$ is of an order not larger than $n^\frac12$.
The upper bound holds as long as $k_n\le cn$ for some $c<1$.)
We also have
$$
k_n^j\le\frac{(2k_n)!}{(2k_n-j)!}\le (2k_n)^j,
\ j=1,...,k_n.
\tag3.8
$$
From (3.6)-(3.8) we have
$$
B(n,k_n,j)\frac{C_1^jj^\frac12}{\Gamma(\frac {j+1}2)}\le
b_2 n^{-j}  (2k_n)^{2j}\ (2k_n)^\frac j2\ \frac{C_1^jj^\frac12}{\Gamma(\frac{ j+1}2)}
\le \frac{j^\frac12C^j}{\Gamma(\frac {j+1}2)}(n^{-1}k_n^\frac52)^j,\tag3.9
$$
for some $C>0$.
Since $\sum_{j=1}^\infty\frac{j^\frac12C^j}{\Gamma(\frac{ j+1}2)}<\infty$,
it follows from (3.9) that
$$
\sum_{j=1}^{k_n}B(n,k_n,j)\frac{C_1^jj^\frac12}{\Gamma(\frac {j+1}2)}
$$
is $O(\frac{k_n^\frac52}n)$ as $n\to\infty$, if $k_n$ is as in part (i), and
is bounded if $k_n$ is as in part (ii).
This proves (3.5-a) and  the part of   (3.5-b) concerning boundedness from $\infty$.

We now turn to (3.5-c) and the part of (3.5-b) concerning boundedness
from 0.
The term in (3.5-b,c)  corresponding to $j=1$
is
$\frac{(k_n!)^2}{\binom n{k_n}^2}
\binom n{2k_n-1}\frac1{(2k_n-1)!}A(k_n-1,1)$.
Define $C(n,k_n)=\frac{(k_n!)^2}{\binom n{k_n}^2}
\binom n{2k_n-1}\frac1{(2k_n-1)!} (2k_n-2)^\frac12\binom {2k_n-2}{k_n-1}^2$.
Using the bound on  $A(k_n-1,1)$
from  Lemma 3, it follows that
for the part of
(3.5-b) concerning boundedness  from 0,
it is enough to show that
$\liminf_{n\to\infty}C(n,k_n)>0$, when $k_n$ is as in part (ii),
and for (3.5-c) it is enough to show that
$\lim_{n\to\infty}C(n,k_n)=\infty$, when $k_n$ is as in part (iii).
Simplifying and making some cancellations, we have
$$
C(n,k_n)=\frac{((n-k_n)!)^2}{n!(n-2k_n+1)!}k_n^4(2k_n-1)^{-2}(2k_n-2)^{\frac12}.
$$
Using this with  (3.7) gives
$$
C(n,k_n)\ge b_1n^{-1}k_n^4(2k_n-1)^{-2}(2k_n-2)^{\frac12}.
$$
Thus, the above stated inequalities indeed hold.
\hfill $\square$
\medskip

\noindent\bf 4. Proofs of Lemmas 2 and 3.

\noindent \bf Proof of Lemma 2.\rm\
The first step of the proof is to develop a probabilistic representation for $A(N,j)$.
Fix $j\ge1$ and $N\ge1$.
Consider two rows each containing
$2N$ spaces.
Randomly fill each of the two rows with $N$ blue balls
and $N$ white balls. Define $X_0=0$, and then
for $m=1,2,...,2N$, use the balls in the first row to
define
$X_m$ as  the number of blue balls in the first $m$ spaces
minus the number of white balls in the first $m$ spaces.
Define $Y_m$ the same way using the balls in the second row.
Then $\{X_m\}$ and $\{Y_m\}$ are independent, and as is well known, each one
has the distribution of  the simple, symmetric one-dimensional random walk,
conditioned to
return to 0 at the $2N$-th step.
Let $U_m=\frac{X_m+Y_m}2$ and  $V_m=\frac{X_m-Y_m}2$.
Then $(U_m,V_m)$ has the distribution of the standard, simple, symmetric two-dimensional
random walk (jumping one unit in each of the four
possible directions with probability $\frac14$), starting from the origin and conditioned to return to the origin
at the $2N$-th step.
To see this, let $\{\Cal X_m\}$ and $\{\Cal Y_m\}$ be independent
copies of  the unconditioned,  simple,
symmetric one-dimensional random walk starting from the origin, and let
$\Cal U_m=\frac{\Cal X_m+\Cal Y_m}2$ and  $\Cal V_m=\frac{\Cal X_m-\Cal Y_m}2$.
Then clearly, $\{\Cal U_m, \Cal V_m\}$ is the unconditioned, simple, symmetric
two-dimensional random walk starting from the origin.
Now $\{U_m,V_m\}$ is equal to $\{\Cal U_m,\Cal V_m\}$ conditioned on
$\Cal X_{2N}=\Cal Y_{2N}=0$, or equivalently, conditioned
on $\Cal U_{2N}=\Cal V_{2N}=0$.

The total number of possible ways of placing $N$ blue balls and $N$
white balls in the first row, and the same number of such balls in the second row
is $\binom{2N}N^2$.
For the moment, fix a set $\{s_r\}_{r=0}^j$ of $j+1$ nonnegative integers
satisfying $\sum_{r=0}^js_r=2N$. Let $t_r=\sum_{i=0}^rs_i$.
Let $D_{s_0,s_1,...,s_j}$ denote the event $\{V_{t_0}=V_{t_1}=...=
V_{t_j}=0\}=\{X_{t_0}=Y_{t_0}, X_{t_1}=Y_{t_1},...,X_{t_j}=Y_{t_j}\}$.
Now for any sequence $\{l_r\}_{r=0}^j$ satisfying $l_r\le s_r$ and
$\sum_{r=0}^jl_r=N$, the probability of the event
$\{X_{t_r}=\sum_{i=0}^r(l_i-(s_i-l_i)),\
r=0,1,...,j\}=
\{X_{t_r}=\sum_{i=0}^r(2l_i-s_i),\
r=0,1,...,j\}$
is $\binom{2N}N^{-1}\prod_{r=0}^j\binom{s_r}{l_r}$,
and thus the probability of the event
$\{X_{t_r}=Y_{t_r}=\sum_{i=0}^r(2l_i-s_i),\
r=0,1,...,j\}$ is
$\binom{2N}N^{-2}(\prod_{r=0}^j\binom{s_r}{l_r})^2$.
Summing now over all possible $\{l_r\}_{r=0}^j$ as above,
it   follows
that
$$
P(D_{s_0,s_1,...,s_j})=
\frac1{\binom{2N}N^2}
\sum\Sb\sum_{r=0}^jl_r=N\\ l_r\le s_r\endSb\
\prod_{r=0}^j \binom{s_r}{l_r}^2.\tag4.1
$$
Letting $m_r=s_r-l_r\ge0$, one sees that
the term involving the summation on the right hand side of (4.1)
can be written as $\sum\prod_{r=0}^j \left(\frac{(l_r+m_r)!}{l_r!m_r!}\right)^2$,
where the sum is over all $\{l_r\}_{r=0}^j$ and $\{m_r\}_{r=0}^j$
satisfying $\sum_{r=0}^jl_r=\sum_{r=0}^jm_r=N$,
and $l_r+m_r=s_r$, for  $r=0,1,...,j$.
Thus, summing (4.1) over all the possible choices of $\{s_r\}_{r=0}^j$, we
obtain
$$
\sum_{\sum_{r=0}^js_r=2N}P(D_{s_0,s_1,...,s_j})=
\frac{A(N,j)}{\binom{2N}N^2}.\tag4.2
$$

\medskip
The next step of the proof is
to estimate
$P(D_{s_0,s_1,...,s_j})$. For this we will need several lemmas.
\medskip

\noindent \bf Lemma 4.\it\ Let $\{Z_n\}_{n=0}^\infty$ be a one-dimensional random walk
which takes jumps of $\pm1$ with probability $\frac14$ each, and remains in its place with
probability $\frac12$. Then
there exit constants $C_1,C_2>0$ such that
$$
\frac{C_1}{n^\frac12}\le P(Z_n=0|Z_0=0)\le \frac{C_2}{n^\frac12},\ \text{for}\ n\ge1.
$$
\medskip
\noindent \bf Proof.\rm\
A direct calculation gives
$$
P(Z_n=0|Z_0=0)=\sum_{i=0}^{[\frac n2]}(\frac12)^{n-2i}(\frac14)^{2i}\binom ni
\binom {n-i}i.
$$
We rewrite this as
$$
P(Z_n=0|Z_0=0)=
(\frac12)^n\sum_{i=0}^{[\frac n2]}(\frac12)^{2i}\binom {2i}i\frac
{\binom ni\binom {n-i}i}{\binom {2i}i}=
(\frac12)^n\sum_{i=0}^{[\frac n2]}(\frac12)^{2i}\binom {2i}i\binom n{2i}.\tag4.3
$$
By Stirling's approximation, there exist positive constants $c_1,c_2$
such that
$$
\frac{c_1}{\sqrt {i+1}}\le(\frac12)^{2i}\binom{2i}i\le \frac{c_2}{\sqrt {i+1}},
\ i=0,1,....\tag4.4
$$
Thus, from (4.3) and (4.4)
we have
$$
c_1\sum_{i=0}^{[\frac n2]}\frac1{\sqrt {i+1}}\binom n{2i}(\frac12)^n\le
P(Z_n=0|Z_0=0)
\le c_2\sum_{i=0}^{[\frac n2]}\frac1{\sqrt {i+1}}\binom n{2i}(\frac12)^n.
\tag4.5
$$
Now let $S_n$ be a random variable distributed according to Binom$(n,\frac12)$. Then we have
$$
\sum_{i=0}^{[\frac n2]}\frac1{\sqrt {i+1}}\binom n{2i}(\frac12)^n=
E(\frac12S_n+1)^{-\frac12}\ 1_{\{S_n \ \text{is even}\}}.\tag4.6
$$
By standard large deviations estimates, $P(|\frac{S_n}n-\frac12|>\epsilon)$
decays exponentially in $n$ for each $\epsilon>0$. Using this along with
(4.5) and (4.6) and leaving to the reader the   little argument to accommodate the
requirement in (4.6) that $S_n$ be even, we conclude that there exist constants
$C_1,C_2>0$ such that
$$
\frac{C_1}{\sqrt {n+1}}\le P(Z_n=0|Z_0=0)\le \frac{C_2}{\sqrt {n+1}}.
$$
\hfill $\square$
\medskip

\noindent \bf Lemma 5.\it\ Let $\{\hat Z_n\}_{n=0}^\infty$ be a simple, symmetric one-dimensional
random walk.

\noindent i. There exists a constant $C_0>0$ such that
$$
P(\hat Z_n=0|\hat Z_0=a)\le \frac {C_0}{\sqrt n}\exp(-\frac{a^2}{2n}),\ \text{for all}\ a\in Z
\ \text{and all}\ n\ge1.
$$

\noindent ii.
Let $L>0$. There exists a constant $c_L>0$ such that for all sufficiently large $n$,
$$
P(\hat Z_{2n}=0|\hat Z_0=2a)\ge \frac{c_L}{\sqrt n}\exp(-\frac{a^2}n), \ \text{for all}\
a\in Z \ \text{satisfying}\ |a|\le Ln^\frac12.
$$
\medskip

\bf \noindent Proof.\rm\ The lemma follows from the local central limit theorem.
It can be proved via a direct calculation, using Stirling's approximation.
(See, for example, \cite{4, page 65}.)
\hfill$\square$
\medskip

\noindent \bf Lemma 6.\it\
Let $\{\hat X_n,\hat Y_n\}_{n=0}^\infty$ be a  simple, symmetric two-dimensional
random walk.

\noindent i. There exist constants $c_1,c_2>0$ such that
$$
 P((\hat X_n,\hat Y_n)=(0,0)|(\hat X_0,\hat Y_0)=(a,0))\le
\frac{c_1}n\exp(-\frac{c_2a^2}n), \ \text{for all}\ a\in Z\ \text{and all}\ n\ge1.
$$
\medskip

\noindent ii. There exists a constant $c_3>0$ such that
$$
P((\hat X_{2n},\hat Y_{2n})=(0,0)|(\hat X_0,\hat Y_0)=(0,0))\ge \frac{c_3}n, \
\text{for all}\ n\ge1.
$$
\noindent \bf Proof.\rm\
Let $H_n$ and $V_n$ denote respectively the number of horizontal and the number of vertical
steps made by the random walk $\{(\hat X_\cdot,\hat Y_\cdot)\}$
during its first $n$ steps.
Then we have
$$
\aligned&P((\hat X_n,\hat Y_n)=(0,0)|(\hat X_0,\hat Y_0)=(a,0))=\\
&\sum_{j+k=n}P(\hat Z_j=0|\hat Z_0=a)P(\hat Z_k=0|\hat Z_0=0)\times P(H_n=j,V_n=k),
\endaligned\tag4.7
$$
where $\{\hat Z_n\}$ is as in Lemma 5.
Since $H_n$ and $V_n$ are each distributed like Binom$(n,\frac12)$, a
standard large deviations estimate gives
$$
P(H_n\ge\frac14n , V_n\ge\frac14n)\ge1-\frac1C\exp(-Cn),\tag4.8
$$
for some $C>0$.
Since $\frac 1{\sqrt j}\exp(-\frac{a^2}{2j})\le \frac2{\sqrt n}\exp(-\frac{a^2}{2n })$,
for $\frac14n\le j\le n$, it follows from (4.7), (4.8) and Lemma 5-i that
$$
 P((\hat X_n,\hat Y_n)=(0,0)|(\hat X_0,\hat Y_0)=(a,0))
 \le
\frac{4C_0^2}n\exp(-\frac{a^2}n)+\frac1C\exp(-Cn).
\tag4.9
$$
Choosing $c_1$ sufficiently large and $c_2>0$ sufficiently small,
part (i)  follows from (4.9) along with the fact that we need only consider $|a|\le n$.

For part (ii), note that
$$
\aligned&P((\hat X_{2n},\hat Y_{2n})=(0,0)|(\hat X_0,\hat Y_0)=(0,0))=\\
&\sum_{2j+2k=2n}P(\hat Z_{2j}=0|\hat Z_0=0)P(\hat Z_{2k}=0|\hat Z_0=0)\times P(H_n=2j,V_n=2k).
\endaligned
$$
Also, we have $P(H_n \ \text{and}\ V_n\ \text{ are even}, H_n\ge\frac14n,V_n\ge\frac14n)\ge C$, for some $C>0$ independent of $n$.
Finally, $P(\hat Z_{2j}=0|\hat Z_0=0)$ can be bounded from below as in
Lemma 5-ii. Part (ii) follows from these observations.
\hfill $\square$

\medskip

We can now estimate $P(D_{s_0,s_1,...,s_j})$.
\medskip

\bf\noindent Lemma 7.\it\
Let $\hat s_r=s_r+1$.
There exists a  constant
$c>0$ such that
$$
P(D_{s_0,s_1,...,s_j})\le(2N+1)^{\frac12}
c^{j+1}(\hat s_0\hat s_1...\hat s_j)^{-\frac12}.\tag4.10
$$
\medskip
\noindent \bf Proof.\rm\  Let $\{\hat X_n, \hat Y_n\}_{n=0}^\infty$
be a simple, symmetric two-dimensional random walk.
Recalling that $t_j=2N$, it follows
by definition that
$$
P(D_{s_0,s_1,...,s_j})=P(\hat Y_{t_0}=\hat Y_{t_1}=\cdot\cdot\cdot=\hat Y_{t_{j-1}}=0|\hat X_0=\hat Y_0=\hat X_{2N}=\hat Y_{2N}=0).
\tag4.11
$$
By the Markov property, we have
$$
\aligned&
P(\hat Y_{t_0}=\hat Y_{t_1}=\cdot\cdot\cdot=\hat Y_{t_{j-1}}=0|\hat X_0=\hat Y_0=\hat X_{2N}=\hat Y_{2N}=0)=\\
&P(\hat Y_{t_0}=0|\hat X_0=\hat Y_0=0)\cdot
P(\hat Y_{t_1}=0|\hat Y_{t_0}=\hat X_0=\hat Y_0=0)\times\cdot\cdot\cdot\times\\
&P(\hat Y_{t_{j-1}}=0|\hat Y_{t_0}=\cdot\cdot\cdot=\hat Y_{t_{j-2}}=\hat X_0=\hat Y_0=0)\times\\
&\frac{ P(\hat X_{2N}=\hat Y_{2N}=0|\hat Y_{t_0}=\cdot\cdot\cdot\hat Y_{t_{j-1}}=\hat X_0=\hat Y_0=0)}
{P(\hat X_{2N}=\hat Y_{2N}=0|\hat X_0=\hat Y_0=0)}.
\endaligned\tag4.12
$$
Note that the process $\{\hat Y_n\}$ in isolation is a one-dimensional random walk distributed
according to the distribution of $\{Z_n\}$ in Lemma 4.
Thus, letting $t_{-1}=0$, we have  from (4.12)
$$
\aligned&
P(\hat Y_{t_0}=\hat Y_{t_1}=\cdot\cdot\cdot=\hat Y_{t_{j-1}}=0|\hat X_0=\hat Y_0=\hat X_{2N}=\hat Y_{2N}=0)=\\
&\frac{ P(\hat X_{2N}=\hat Y_{2N}=0|\hat Y_{t_0}=\cdot\cdot\cdot\hat Y_{t_{j-1}}=\hat X_0=\hat Y_0=0)}
{P(\hat X_{2N}=\hat Y_{2N}=0|\hat X_0=\hat Y_0=0)}\times
\prod_{k=0}^{j-1}P(Z_{t_k}=0|Z_{t_{k-1}}=0).
\endaligned\tag4.13
$$
Recall that $t_k-t_{k-1}=s_k$ and $s_j=2N-t_{j-1}$. Let $s_k'=s_k$,
if $s_k\ge1$, and $s_k'=1$, if $s_k=0$. Since
$P(Z_{t_k}=0|Z_{t_{k-1}}=0)=P(Z_{t_k-t_{k-1}}=0|Z_0=0)$ it follows
from Lemma 4 that $P(Z_{t_k}=0|Z_{t_{k-1}}=0)\le
\frac{C_2}{\sqrt{s_k'}}$. From Lemma 6-ii,
it follows that $P(\hat X_{2N}=\hat Y_{2N}=0|\hat X_0=\hat
Y_0=0)\ge \frac{c_3}N$. Using these facts along with (4.11) and
(4.13), it follows that  if we show that
$$
P(\hat X_{2N}=\hat Y_{2N}=0|\hat Y_{t_0}=\cdot\cdot\cdot\hat Y_{t_{j-1}}=\hat X_0=\hat Y_0=0)\le \frac
{C^j}{N^\frac12(2N-t_{j-1})^\frac12},\tag4.14
$$
for some $C>0$ and $s_j=2N-t_{j-1}\ge1$, then we will obtain (4.10)
with $\hat s_r$ replaced by $s_r'$. By increasing the constant $c$
in (4.10), one can always replace $s_r'$ by $\hat s_r$. Thus, it
remains to prove (4.14).

By the Markov property, we have
$$
\aligned&
P(\hat X_{2N}=\hat Y_{2N}=0|\hat Y_{t_0}=\cdot\cdot\cdot\hat Y_{t_{j-1}}=\hat X_0=\hat Y_0=0)=\\
&\sum_{a\in R} P(\hat X_{2N-t_{j-1}}=\hat Y_{2N-t_{j-1}}=0|\hat X_0=a,\hat Y_0=0)\times
\nu(a),\endaligned\tag4.15
$$
where
$$\nu(a)=P(\hat X_{t_{j-1}}=a|\hat
Y_{t_0}=\cdot\cdot\cdot\hat Y_{t_{j-1}}=\hat X_0=\hat Y_0=0).
$$
As in the proof of Lemma 6,
let $H_n$ denote the number of horizontal steps taken by the random walk $\{\hat X(\cdot),\hat Y(\cdot)\}$
during its first $n$ steps.
Let
$$
\mu(m)=P(H_{t_{j-1}}=m|\hat
Y_{t_0}=\cdot\cdot\cdot\hat Y_{t_{j-1}}=\hat X_0=\hat Y_0=0),
$$
and let $W$ be distributed like $\mu$.
Then $\nu$ is distributed like $\hat Z_W$, where $\{\hat Z_n\}$ is a simple, symmetric one-dimensional random walk,
starting from 0 and independent of $W$.
We will show later that for some $\gamma,C>0$,
$$
\aligned&\mu([0,\gamma t_{j-1}])=P(H_{t_{j-1}}\le\gamma t_{j-1}|\hat
Y_{t_0}=\cdot\cdot\cdot\hat Y_{t_{j-1}}=\hat X_0=\hat Y_0=0)=\\
&P(W\le \gamma t_{j-1})\le \frac {C^j}\gamma\exp(-\gamma t_{j-1}).
\endaligned
\tag4.16
$$

By Lemma 5-i, it follows that
$$
P(\hat Z_n=a|\hat Z_0=0)\le \frac{C_0}{\sqrt{\gamma t_{j-1}}}\exp(-\frac{a^2}{2t_{j-1}}),
\ \text{for}\  \gamma t_{j-1}\le n\le t_{j-1}.\tag4.17
$$
From (4.16) and (4.17) we conclude that
$$
\nu(a)=P(\hat Z_W=a)\le
\frac{C_0}{\sqrt{\gamma t_{j-1}}}\exp(-\frac{a^2}{2t_{j-1}})+\frac{C^j}\gamma\exp(-\gamma t_{j-1}).\tag4.18
$$
Since $\nu(a)=0$, if $a>t_{j-1}$, it follows from (4.18) that
$$
\nu(a)\le \frac{k_1^j}{\sqrt{t_{j-1}}}\exp(-\frac{k_2a^2}{t_{j-1}}),\tag4.19
$$
for some $k_1,k_2>0$.
From (4.15), (4.19) and Lemma 6-i, we obtain
$$
\aligned& P(\hat X_{2N}=\hat Y_{2N}=0|\hat Y_{t_0}=\cdot\cdot\cdot\hat Y_{t_{j-1}}=\hat X_0=\hat Y_0=0)
\le\\
& \sum_{a\in R}\frac{c_1}{2N-t_{j-1}}\exp(-\frac{c_2a^2}{2N-t_{j-1}})\frac{k^j_1}{\sqrt{t_{j-1}}}
\exp(-\frac{k_2a^2}{t_{j-1}}).
\endaligned\tag4.20
$$
For an appropriate $\hat C>0$, the right hand side of (4.20) can be
bounded from above by $\hat C^j\int_{-\infty}^\infty
\frac1{(2N-t_{j-1})\sqrt{t_{j-1}}}\exp(-\frac{c_2x^2}{2N-t_{j-1}})\exp(-\frac{k_2x^2}{t_{j-1}})dx$.
Evaluating this integral gives the estimate in (4.14). Thus, to
complete the proof of the lemma, it remains to prove (4.16).

We mention that it is intuitive that
$P(H_{t_{j-1}}\le\gamma t_{j-1}|\hat
Y_{t_0}=\cdot\cdot\cdot\hat Y_{t_{j-1}}=\hat X_0=\hat Y_0=0)\le P(H_{t_{j-1}}\le \gamma t_{j-1})$,
and this would then give (4.16).
The intuition comes from the fact that the smaller $H_{t_{j-1}}$ is, the more moves
$\{\hat Y_n\}$ makes, and the more moves $\{\hat Y_n\}$ makes, the more difficult it is for
it to
have the prescribed zeroes. However, a proof of this is rather complicated and quite tedious.
It turns out that a rather crude estimate will suffice in order to obtain (4.16).
We have
$$
\aligned&\mu([0,\gamma t_{j-1}])=
P(H_{t_{j-1}}\le\gamma t_{j-1}|\hat
Y_{t_0}=\cdot\cdot\cdot\hat Y_{t_{j-1}}=\hat X_0=\hat Y_0=0)=\\
&\frac{P(H_{t_{j-1}}\le\gamma t_{j-1}, \hat
Y_{t_0}=\cdot\cdot\cdot\hat Y_{t_{j-1}}=\hat X_0=\hat Y_0=0)}
{P(\hat
Y_{t_0}=\cdot\cdot\cdot\hat Y_{t_{j-1}}=\hat X_0=\hat Y_0=0)}\le\frac{P(H_{t_{j-1}}\le\gamma t_{j-1})}{P(\hat
Y_{t_0}=\cdot\cdot\cdot\hat Y_{t_{j-1}}=0|\hat Y_0=0)}.
\endaligned\tag4.21
$$
By a standard large deviations estimate,
$$
P(H_{t_{j-1}}\le\gamma t_{j-1})\le c\exp(-l_\gamma t_{j-1}), \ \text{where}\
\lim_{\gamma\to0}l_\gamma=\log 2.\tag4.22
$$
(To see this, note that $P(H_{t_{j-1}}=0)=(\frac12)^{t_{j-1}}=\exp(-(\log 2) t_{j-1})$.)
By Lemma 4, we have
$$
P(\hat Y_{t_0}=\cdot\cdot\cdot\hat Y_{t_{j-1}}=0|\hat Y_0=0)
=\prod_{k=0}^{j-1}P(\hat Y_{t_k}=0|\hat Y_{t_{k-1}}=0)\ge
\frac{C_1^j}{(s_0's_1'\cdot\cdot\cdot s_{j-1}')^\frac12}.
$$
Since the $\{s_k'\}$ satisfy $\sum_{k=0}^{j-1}s_k'\le t_{j-1}+j$, it follows that
$\sup_{\{s_k'\}}s_0's_1'\cdot\cdot\cdot s_{j-1}'\le (\frac{t_{j-1}+j}j)^j\le\exp(t_{j-1})$.
Thus,
$$
P(\hat Y_{t_0}=\cdot\cdot\cdot\hat Y_{t_{j-1}}=0|\hat Y_0=0)
\ge C_1^j\exp(-\frac12 t_{j-1}).\tag4.23
$$
Now (4.16) follows from (4.21)-(4.23) along with the fact that $\log 2>\frac12$.
\hfill $\square$

We can now complete the proof of Lemma 2. From (4.2) and (4.10), we have
$$
\frac{A(N,j)}{\binom{2N}N^2}\le
(2N+1)^{\frac12}
c^{j+1}\sum_{\sum_{r=0}^js_r=2N}
(\hat s_0\hat s_1...\hat s_j)^{-\frac12}.
$$
Let $\hat S_j=\hat s_j+j-1$. Then $(\hat s_j)^{-\frac12}=
(\hat S_j)^{-\frac12}(\frac{\hat S_j}{\hat
s_j})^\frac12\le\frac {j^\frac12}{(\hat S_j)^\frac12}$.
Thus it follows from the above inequality that
$$
\frac{A(N,j)}{\binom{2N}N^2}\le
(2N+1)^{\frac12}
c^{j+1}j^\frac12\sum_{\sum_{r=0}^js_r=2N}
(\hat s_0\hat s_1...\hat s_{j-1}\hat S_j)^{-\frac12}.\tag4.24
$$
The replacement of $\hat s_j$ by $\hat S_j$ was made for technical
reasons which will become clear below.
Making the substitutions $x_r=\frac{s_r}{2N}$ and  $\hat x_r
=\frac{\hat s_r}{2N}$, for $r=0,...,j$, and $\hat X_j=\frac{\hat S_j}{2N}$,
we rewrite the right hand
side of (4.24) as
$$
(2N+1)^{\frac12}
c^{j+1}j^\frac12(2N)^{\frac{j-1}2} \sum\Sb\sum_{r=0}^jx_r=1\\ (2N)x_r\ \text{is a nonnegative integer}\endSb
(\hat x_0\hat x_1...\hat x_{j-1}\hat X_j)^{-\frac12}(2N)^{-j}.\tag4.25
$$

Let $C_{x_0,x_1,...,x_{j-1}}$ denote the hyper-cube
$\prod_{r=0}^{j-1}[x_r,x_r+\frac1{2N}]=\prod_{r=0}^{j-1}[x_r,\hat x_r]$. Consider $\cup
C_{x_0,x_1,...,x_{j-1}}$, where the union is over
all $\{x_0,...,x_{j-1}\}$ for which $(2N)x_r$ is a nonnegative integer and
$\sum_{r=0}^{j-1}x_r\le 1$. This union
is contained in $V_{1+\frac j{2N}}
\equiv\{(y_0,y_1,...,y_{j-1}):y_r\ge0,\ \sum_{r=0}^{j-1}y_r\le 1+\frac j{2N}\}$.
We have
$$
\aligned&(\hat x_0\hat x_1...\hat x_{j-1}\hat
X_j)^{-\frac12}\le (y_0y_1...y_{j-1}y_j)^{-\frac12},\ \text{for
all}\ (y_0,y_1,...,y_{j-1})\in C_{x_0,x_1,...,x_{j-1}},\\
& \text{where}\ y_j=1+\frac
j{2N}-y_0-y_1-...-y_{j-1}.\endaligned\tag4.26
$$
To see that (4.26) holds, note that $\hat x_r\ge
y_r, \ r=0,...,j-1$, for $(y_0,y_1,...,y_{j-1})\in
C_{x_0,x_1,...,x_{j-1}}$.  Also,
$$
\aligned&\hat X_j=\frac{\hat S_j}{2N}=\frac{\hat
s_j+j-1}{2N}=\frac{s_j+j}{2N}
=\frac{2N+j-s_0-s_1-...-s_{j-1}}{2N}\\
&=1+\frac j{2N}-x_0-x_1-...-x_{j-1}\ge
1+\frac j{2N}-y_0-y_1-...-y_{j-1},\\
& \text{for}\ (y_0,y_1,...,y_{j-1})\in C_{x_0,x_1,...,x_{j-1}}.\endaligned
$$

In light of these facts it follows that the sum on the right hand
side of (4.25) is dominated by a certain lower Riemann sum for
 $\int_{\Cal S_{1+\frac j{2N}}}(\prod_{r=0}^jy_r)^{-\frac12}dy_0dy_1...
dy_{j-1}$, where
$\Cal S_\lambda=\{(y_0,y_1,...,y_{j-1}):y_r\ge0,
\sum_{r=0}^{j-1}y_r\le\lambda\}$.
Replacing the sum in (4.25) with this integral, and substituting
the resulting expression  into the right hand
side of (4.24) gives
$$
\frac{A(N,j)}{\binom{2N}N^2}\le
(2N+1)^{\frac12}
c^{j+1}j^\frac12(2N)^{\frac{j-1}2}
\int_{\Cal S_{1+\frac j{2N}}}(\prod_{r=0}^jy_r)^{-\frac12}dy_0dy_1...
dy_{j-1}.\tag4.27
$$
A change of variables shows that
$$
\int_{\Cal S_{1+\frac j{2N}}}(\prod_{r=0}^jy_r)^{-\frac12}dy_0dy_1...
dy_{j-1}=(1+\frac j{2N})^{\frac{j-1}2}
\int_{\Cal S_1
}(\prod_{r=0}^jy_r)^{-\frac12}dy_0dy_1...
dy_{j-1}.\tag4.28
$$
As is well-known from the theory of Dirichlet distributions,
$$
\int_{\Cal S_1}(\prod_{r=0}^jy_r)^{-\frac12}dy_0dy_1...
dy_{j-1}=\frac{\pi^{\frac {j+1}2}}{\Gamma(\frac {j+1}2)}.\tag4.29
$$
From (4.27)-(4.29) we conclude that
$$
\frac{A(N,j)}{\binom{2N}N^2}\le (2N+1)^{\frac12}
c^{j+1}j^\frac12(2N)^{\frac{j-1}2} (1+\frac j{2N})^{\frac{j-1}2}
\frac{\pi^{\frac {j+1}2}}{\Gamma(\frac {j+1}2)}.\tag4.30
$$
The inequality for $A(N,j)$ in Lemma 2 follows from (4.30).
It is trivial to check that  $A(N,j)$ is increasing in $N$;
thus, the inequality for $A(k-j,j)$ in Lemma 2 holds as stated.
\hfill $\square$

\medskip
\noindent\bf Proof of Lemma 3.\rm\
To prove the lemma we will need the following lemma, which complements Lemma 6.

\medskip

\noindent\bf Lemma 8.\it\
Let $\{\hat X_n,\hat Y_n\}_{n=0}^\infty$ be a  simple, symmetric two-dimensional
random walk. Let $L>0$. There exist constants $c_{L,1},c_{L,2}>0$ such that for all sufficiently
large $n$
$$
\aligned& P((\hat X_{2n},\hat Y_{2n})=(0,0)|(\hat X_0,\hat Y_0)=(2a,0))\ge
\frac{c_{L,1}}n\exp(-\frac{c_{L,2}a^2}n), \\
&\ \text{for }\ a\in Z\ \text{satisfying}\ |a|\le L\sqrt n.
\endaligned
$$
\medskip

\noindent \bf Proof.\rm\ Let $H_n,V_n$ be as in the proof of Lemma
6, and let $\{\hat Z_n\}$ be as in Lemma 5. We have
$$
\aligned&P((\hat X_{2n},\hat Y_{2n})=(0,0)|(\hat X_0,\hat Y_0)=(2a,0))=\\
&\sum_{j+k=n}P(\hat Z_{2j}=0|\hat Z_0=2a)P(\hat Z_{2k}=0|\hat
Z_0=0)\times P(H_n=j,V_n=k).
\endaligned\tag4.31
$$
The proof of the lemma follows easily from (4.31), (4.8) and Lemma 5-ii.
\hfill $\square$
\medskip

We can now prove Lemma 3.
Let $\{(\hat X_n,\hat Y_n)\}$ denote a
simple, symmetric two-dimensional random walk.
Using the notation in the proof of Lemma 2, but with $k-1$ in place of $N$,
recall that for $j=1$,
$$
P(D_{s_0,s_1})=P(\hat Y_{s_0}=0|\hat X_0=\hat Y_0=\hat X_{2k-2}=\hat Y_{2k-2}=0),
$$
and thus from (4.2),
$$
\frac{A(k-1,1)}{\binom{2k-2}{k-1}^2}=
\sum_{l=0}^{2k-2}P(\hat Y_l=0|\hat X_0=\hat Y_0=\hat X_{2k-2}=\hat Y_{2k-2}=0).\tag4.32
$$
For $m$ satisfying $[\frac14 k]\le m\le [\frac34 k]$, we have
$$
\aligned&P(\hat Y_{2m}=0|\hat X_0=\hat Y_0=\hat X_{2k-2}=\hat Y_{2k-2}=0)\ge\\
&\sum_{r=-[\sqrt k]}^{[\sqrt k]}P(\hat Y_{2m}=0,\hat X_{2m}=2r|\hat X_0=\hat Y_0=0)\times\\
&\frac{P(\hat X_{2k-2}=
\hat Y_{2k-2}=0|\hat X_0=\hat Y_0=\hat Y_{2m}=0,\hat X_{2m}=2r)}
{P(\hat X_{2k-2}=\hat Y_{2k-2}=0|\hat X_0=\hat Y_0=0)}.\endaligned\tag4.33
$$
We have
$P(\hat X_{2k-2}=
\hat Y_{2k-2}=0|\hat X_0=\hat Y_0=\hat Y_{2m}=0,\hat X_{2m}=2r)
=P(\hat X_{2k-2-2m}=\hat Y_{2k-2-2m}=0|\hat X_0=2r,\hat Y_0=0)$.
Thus, in light of the above-specified range of $m$ and of $r$, it follows from Lemma 8
and Lemma 6-i that for sufficiently large $k$,
$\frac{P(\hat X_{2k-2}=\hat Y_{2k-2}=0|\hat X_0=\hat Y_0=\hat Y_{2m}=0,\hat X_{2m}=2r)}
{P(\hat X_{2k-2}=\hat Y_{2k-2}=0|\hat X_0=\hat Y_0=0)}$ is bounded from below by a positive
constant. By Lemma 8  and the above-specified bound on $m$,  it also follows that
$P(\hat Y_{2m}=0,\hat X_{2m}=2r|\hat X_0=\hat Y_0=0)$ is bounded from below by $\frac Ck$,
for some $C>0$.
Thus, we conclude from (4.33) that for sufficiently large $k$,
$$
\aligned&P(\hat Y_{2m}=0|\hat X_0=\hat Y_0=\hat X_{2k-2}=\hat Y_{2k-2}=0)\ge C_1k^{-\frac12},\\
& \text{ for
some}\ C_1>0\ \text{and for }\ m\ \text{satisfying}\ [\frac14k]\le m\le [\frac34 k].
\endaligned\tag4.34
$$
Now (4.32) and (4.34) give
$$
\frac{A(k-1,1)}{\binom{2k-2}{k-1}^2}\ge
C_2k^\frac12,
$$
for
some  $C_2>0$ and $k$ sufficiently large.
This is clearly equivalent to the lemma.
\hfill $\square$

\medskip

\noindent \bf Acknowledgement.\rm\ The author thanks the referees for their careful reading of the
paper. He  thanks one of them in particular for pointing out an error in a previous
version of the paper.

\end